\UseRawInputEncoding 
\documentclass[review,3p,10pt]{elsarticle}
\usepackage{amsmath}
\usepackage[mathlines]{lineno}
\usepackage{amssymb}
\usepackage{amsmath}
\usepackage{cases}
\usepackage{epstopdf}
\usepackage[caption=false,font=footnotesize]{subfig}
\usepackage{fixltx2e}
\usepackage{amssymb}
\usepackage{mathtools}
\usepackage{bm}
\usepackage{multirow}
\usepackage{color}
\usepackage{algorithm}
\usepackage[noend]{algpseudocode}
\usepackage{algorithmicx,algorithm}
\usepackage{amsthm}
\usepackage{mathrsfs}

\newcommand*\patchAmsMathEnvironmentForLineno[1]{%
  \expandafter\let\csname old#1\expandafter\endcsname\csname #1\endcsname
  \expandafter\let\csname oldend#1\expandafter\endcsname\csname end#1\endcsname
  \renewenvironment{#1}%
     {\linenomath\csname old#1\endcsname}%
     {\csname oldend#1\endcsname\endlinenomath}}%
\newcommand*\patchBothAmsMathEnvironmentsForLineno[1]{%
  \patchAmsMathEnvironmentForLineno{#1}%
  \patchAmsMathEnvironmentForLineno{#1*}}%
\AtBeginDocument{%
\patchBothAmsMathEnvironmentsForLineno{equation}%
\patchBothAmsMathEnvironmentsForLineno{align}%
\patchBothAmsMathEnvironmentsForLineno{flalign}%
\patchBothAmsMathEnvironmentsForLineno{alignat}%
\patchBothAmsMathEnvironmentsForLineno{gather}%
\patchBothAmsMathEnvironmentsForLineno{multline}%
}

\journal{Control Engineering Practice}

\begin{document}
%
\begin{frontmatter}


\title{Deep Learning-Aided Model Predictive Control of Wind Farms for AGC Considering the Dynamic Wake Effect}

\author[label1]{Kaixuan~Chen}
\author[label1]{Jin~Lin\corref{cor1}}
\ead{linjin@tsinghua.edu.cn}
\author[label1]{Yiwei~Qiu}
\author[label1]{Feng Liu}
\author[label1,label3]{Yonghua~Song}

\address[label1]{State Key Laboratory of Control and Simulation of Power Systems and Generation Equipment, Department of Electrical Engineering, Tsinghua University, Beijing 100087, China}
\address[label3]{Department of Electrical and Computer Engineering, University of Macau, Macau 999078, China}

\cortext[cor1]{Corresponding author}
\address{}

\begin{abstract}
To provide automatic generation control (AGC) service, wind farms (WFs) are required to control their operation dynamically to track the time-varying power reference. Wake effects impose significant aerodynamic interactions among turbines, which remarkably influence the WF dynamic power production. 
The nonlinear and high-dimensional nature of dynamic wake model, however, brings extremely high computation complexity and obscure the design of WF controllers. 
This paper overcomes the control difficulty brought by the dynamic wake model by proposing a novel control-oriented reduced order WF model and a deep-learning-aided model predictive control (MPC) method.
Leveraging recent advances in computational fluid dynamics (CFD) to provide high-fidelity data that simulates WF dynamic wake flows, two novel deep neural network (DNN) architectures are specially designed to learn a dynamic WF reduced-order model (ROM) that can capture the dominant flow dynamics. Then, a novel MPC framework is constructed that explicitly incorporates the obtained WF ROM to coordinate different turbines while considering dynamic wake interactions.
The proposed WF ROM and the control method are evaluated in a widely-accepted high-dimensional dynamic WF simulator whose accuracy has been validated by realistic measurement data. A 9-turbine WF case and a larger 25-turbine WF case are studied. By reducing WF model states by many orders of magnitude, the computational burden of the control method is reduced greatly. Besides, through the proposed method, the range of AGC signals that can be tracked by the WF in the dynamic operation is extended compared with the existing greedy controller.  
\end{abstract}


\end{frontmatter}

\section{Introduction}

The increasing penetration of wind energy in the power system is an irreversible trend \cite{veers2019grand}. Instead of just maximizing their production, wind farms (WFs) are required to provide many auxiliary services for the power system in high-renewable era \cite{veers2019grand}. Automatic generation control (AGC), which requires the WF power generation to track a time-varying reference given by power system operators during a time span of several minutes, is one of the main auxiliary services that has already become a requirement for WFs in some jurisdictions\cite{veers2019grand, aho2016active}. 

Active power tracking at a stand-alone turbine level has been investigated thoroughly \cite{aho2016active, 2019Toward}. However, control of dynamic power production at the farm level is much more complicated because of the aerodynamic coupling among turbines that occurs at the time scale of minutes comparable to those of the AGC signals \cite{shapiro2017model}. For example, if an upstream wind turbine (WT) increases its power generation, the overall WF power production may not increase because the downstream turbines experience a reduced wind speed, which is known as the wake effect \cite{boersma2017tutorial, YANG201537}. More importantly, such wake effects have a minutes-level dynamic process because the reduced wind availability will travel from the upstream turbine to downstream turbines through mean flow advection. Such an aerodynamic interaction has an important effect on WF dynamic power production and has been demonstrated by high-fidelity simulation \cite{fleming2016computational}, wind tunnel experiments \cite{campagnolo2016wind} and real field tests \cite{fleming2017field}.

WF control considering the wake effect has attracted much attention recently. Many studies use static parametric wake models to maximize WF power production \cite{article}. Though these studies provide the WF maximal power extraction when operating at the steady-state optimal point, they cannot optimize the dynamic control signals. Such static wake models that ignore the time-varying nature of the wake interactions are inadequate for AGC studies where a dynamic control trajectory is required to control WFs that follow the time-varying power reference and interact with flow dynamics \cite{shapiro2017model, boersma2017tutorial}.

For WF power tracking, the problem lies in how to dynamically adjust the power demand distribution among different turbines \cite{vali2019active}. 
{\color{black}Literature} \cite{sorensen2005wind} proposes the currently widely used proportion distribution (ProD), which distributes the power command proportionally to the available power capabilities of each WT. {\color{black}Literature} \cite{van2017active} and \cite{vali2019active} uses a classical feedback control to improve the power tracking performance. 
However, because no optimization is performed, the range of trackable power is much lower than the available power of the WF when it operates at its global optimal point considering wake interactions  \cite{munters2018dynamic}. To optimize the WF power tracking performance, closed loop control with a dynamical wake model is required.

The dynamical behavior of the wind flow in a WF is governed by the unsteady Navier-Stokes (NS) equations \cite{boersma2017tutorial}. 
Recent advances in computational fluid dynamics (CFD) have enabled high-fidelity dynamic WF modeling based on the numerical discretization of the NS equations \cite{meyers2010large, boersma2018control,deskoswake}. 
A series of optimal WF control studies that apply model predictive control (MPC) techniques directly in the CFD model have been provided by \cite{munters2018dynamic, vali2019adjoint}. In these studies, the NS-equation based dynamic WF model is numerically discretized as state-space representations, and control trajectories are optimized by adjoint-based gradient methods. The range of trackable power is improved to the maximal WF power extraction considering the wake effect\cite{vali2019adjoint}. However, the nonlinearity and extremely high dimensions of the dynamic WF model make it computationally prohibitive for control implementation as stated by the authors themselves \cite{boersma2017tutorial, munters2018dynamic}, which have $10^4$ or more states and take up days of distributed computation. 

Therefore, to consider the dynamic wind flow interactions, a computationally-acceptable optimal control method is not yet available in the current literature and needs to be addressed.

With the advances in CFD, many design and validation processes generate a high volume of CFD data, opening the possibility of using data-driven approaches to facilitate the controller design \cite{kutz2017deep}. Over the past years, deep learning approaches have shown great successes in learning efficient low-order representations from high-dimensional data to reduce control complexity in some research fields such as image-based robot control \cite{watter2015embed}, fluid fast simulations \cite{hennigh2017lat}\cite{zhang2020novel} and building temperature control \cite{terzi2020learning}. Compared to optimal control based on physical models, the data-driven modelling is promising for reducing the model complexity and reducing the computational cost for control applications \cite{zhang2020novel}\cite{GUO2021104758}\cite{HEDINGER201917}. 
For this reason, in this paper we propose a novel, data-driven approach to solve the optimal control problem in the high-dimensional nonlinear WF system.  
Particularly, a control-oriented reduced-order model (ROM) is identified and a novel MPC framework is designed. Two novel deep neural network (DNN) architecture and corresponding loss functions are constructed to learn both WF reduced-order dynamics and power equations in a specially designed form. Then a novel MPC algorithm is proposed that explicitly incorporates the learned WF ROM to coordinate different turbines considering their dynamic wake interactions. The main contributions of this paper are summarized as follows.

\vspace{10pt}
1) A specially-designed autoencoder is constructed to learn WF reduced-order dynamics from CFD data in a globally linear form. Then, a separate fully connected DNN is designed to learn the local linearization of reduced-order power output equations. The novel control-oriented WF reduced-order model is developed in the data-driven approach, which considers both the accuracy of model reduction and the control compatibility. The number of model states is reduced by orders of magnitude. More importantly, the specially designed network allows the learned WF ROM to be readily incorporated into a MPC framework.

2) A novel control framework embedding the derived WF ROM is proposed, through which the WF dynamic power production can be optimized in a computational effective manner. As far as the authors know, this is the first optimal WF controller that optimizes the dynamic wake interactions and the computation burden is practical for control applications. The proposed method is tested in a high-dimensional WF simulator whose accuracy has been validated by realistic measurement data and multiple previous studies. By coordinating different turbines, the range of AGC signals that the WF can track in the dynamic operation is increased from the greedy power to the maximal extractable power considering wake interactions. The computational advantages of the proposed control method are analyzed theoretically and experimentally.

\vspace{10pt}

The remainder of this paper is organized as follows. Section \ref{wakemodel} reviews the dynamic WF model. Section \ref{network} constructs the DNN architecture to learn the WF ROM. Section \ref{control} introduces the deep learning-aided MPC framework. Case studies and conclusions are presented in Section \ref{case} and \ref{cons}.


\section{Dynamic Wind Farm Model} \label{wakemodel}
In this section, we briefly review the dynamic WF model that is commonly employed in dynamic WF simulations \cite{meyers2010large, boersma2018control}. The WF flow field is governed by the incompressible Navier-Stokes equations: 
\begin{align}
   \frac{\partial  \bm{v}  } {\partial t} + (\bm{v} \cdot \nabla) \bm{v}=& - \frac{1}{\rho} \nabla p + \eta \nabla^2 \bm{v}+ \frac{1}{\rho}(\sum_{i=1}^{N} \bm{f}_i)  \label{equ:ns}\\
   \nabla \cdot \bm{v}=&\hspace{3pt}0  \label{equ:conti}
\end{align}
where $\bm{v}=[v_x,v_y]$ represents the velocity field at the hub height with  $v_x, v_y$ denoting the streamwise and spanwise velocity components, respectively;  $p$ is the pressure field. The two constants $\eta$ and $\rho$ denote the kinematic viscosity and the air density, respectively; $N$ is the number of turbines in the WF. The classical actuator disk model is utilized to model the turbines in the dynamic wind flow field, which exerts a thrust force into the wind flow and extracts energy \cite{meyers2010large, boersma2018control}. $\bm{f}_i$ represents the thrust force of turbine $i$, which acts on the position where the turbine rotor disk is located: 
\begin{align}
  \bm{f}_i= - \frac{1}{2} \rho U_i^2 C_{T_i}^{\prime} \bm{e}_{\bot,i} \label{wtforce} 
\end{align}
where $C_{T_i}^\prime$ is the disk-based thrust coefficient, which can directly be related to the blade pitch angle and rotational speed \cite{meyers2010large}; 
$\bm{e}_{\bot,i}$ is the rotor-perpendicular unit vector that defines the orientation of the rotor disk and can be controlled by the turbine yaw angle $\gamma_i$:
\begin{align}
\bm{e}_{\bot,i} = \bm{e}_x \cos(\gamma_i) + \bm{e}_y \sin(\gamma_i).
\end{align}
where $\bm{e}_x$ and $\bm{e}_y$ represent the streamwise and spanwise unit vector, respectively.
$U_i$ is the rotor-perpendicular velocity averaged over the rotor disk. 
As commonly employed in wind farm studies, $C_{T_i}^\prime$ and $\gamma_i$ of each turbine are considered as the control variables and are used to regulate WF operation \cite{meyers2010large, boersma2018control, vali2019adjoint}.

The power generated by a turbine is computed as follows:
\begin{align}
  P^i=\frac{1}{2} \rho A_d  U_i^3 C_{P_i}   \label{equ:power}
\end{align}
where $P^i$ is the produced power of the $ith$ turbine. The power coefficient $C_{P_i} = c_p C_{T_i}^\prime$. The constant $c_p$  is the loss factor and set to 0.9 as discussed in \cite{munters2017optimal}. $A_d$ is the rotor area. 

No analytic solution has been to be found for these partial differential equations (\ref{equ:ns})--(\ref{equ:conti}). The simulation of wind flow is typically solved by CFD.
Standard CFD methods such as large eddy simulation are used to discretize the partial differential equation model (\ref{equ:ns})--(\ref{equ:conti}) spatially and temporally to allow for high-fidelity numerical simulation \cite{boersma2018control}. The wind farm field is spatially discretized over a grid of $(N_x \times N_y)$ cells with $N_x$ points in the streamwise direction and $N_y$ points in the spanwise direction. After spatially discretizing, each point of this grid has its own equations that describe the flow velocity and pressure at that point. Then, the discretized WF model 
can be represented in a high-dimensional nonlinear descriptor state-space form as follows: 
\begin{align}
  \textbf{x}_{t+1}=f(\textbf{x}_{t}, \textbf{u}_{t}) \label{CFDstate} \\
  \textbf{P}_{t}=h(\textbf{x}_{t}, \textbf{u}_{t})
\end{align}
where $t$ denotes the time steps after temporally discretizing. The wind velocity fields $v_x$ and $v_y$ are discretized over the grid. Consequently, the total velocity state number for the dynamic WF after discretization is $(N_x \times N_y  \times 2)$ with 2 velocity components per grid point. To save space, we use $\textbf{x}_t$ to denote the full system states at time step $t$, which is a
 tensor of size $(N_x,N_y, 2)$ and describes the velocity snapshot of the whole WF grid at the time instant $t$. $\textbf{u}_{t}=[C_{T_1}^\prime, \gamma_1, \ldots,C_{T_{N}}^\prime, \gamma_{N}]$ is the control input to the WF. $\textbf{u}\in\mathbb{R}^{2N\times1}$ consist of $[C_{T}^\prime, \gamma]$ of each turbine. $\textbf{P}_{t}=[P^1_{t},\ldots,P^{N}_{t}] \in\mathbb{R}^{N\times1} $ is the power production of each turbine in the WF, which is considered as the output of the system. The relations $f$ and $h$ are nonlinear, whose concrete form depends on the specific discretizing method. We omit the details of the discretization procedure, which is the focus of CFD studies and outside the scope of the current paper. Readers interested in the details of CFD discretization are referred to \cite{boersma2018control, munters2017optimal}.

\section{A Control-Oriented WF ROM Based on Deep Learning} \label{network}
Though the standard dynamic WF model can capture the wake dynamics, the number of states can easily reach $10^{4}$ or more due to the enormous number of spatial states required for simulation, which is excessively computationally expensive for controller design and analysis \cite{boersma2017tutorial, munters2018dynamic}. Ideally, a control-oriented model is desired to capture the dominant flow dynamics in a computationally efficient manner. This is the motivation for the ROM.
In this section, we formulate the ROM problem. Two DNN architectures for learning the WF ROM from the high-dimensional CFD data is proposed, based on which the computationally tractable controller can be designed. 

\subsection{Problem Formulation} \label{probform}
Though the dynamic WF model is extremely high dimensional in state space, it is observed that the dominant dynamics of spatial-temporal flow data always evolve in low-rank subspaces spanned by a few spatial coordinates \cite{kutz2017deep, annoni2017method, lusch2018deep}. Our goal is to infer a dynamic WF ROM in which optimal control can be readily performed. To this end, intrinsic coordinates $\textbf{z}=\phi(\textbf{x})$ that dominate the WF flow data need to be identified so that high-dimensional $\textbf{x}$ can be mapped to low-dimensional vectors $\textbf{z} \in\mathbb{R}^{n_z}$ with $n_z\ll (N_x\times N_y) $. At the same time, the system dynamics with control are also required to be identified in the latent space $\widehat{\textbf{z}}_{t+1} = f^{\textrm{ROM}}(\textbf{z}_t,\textbf{u}_t)$ so that the control problem can be solved in reduced-order space $\textbf{z}$ instead of the original full state space $\textbf{x}$. 

To this end, the ROM has to fulfill three properties: (1) \textit{Reconstruction}. The reduced-order states $\textbf{z}=\phi(\textbf{x})$ must capture sufficient information about the full states $\textbf{x}_t$ (enough to enable reconstruction). (2) \textit{Future state prediction}. The dynamic ROM $f^{\textrm{ROM}}$ must allow for accurate prediction of the next few reduced-order states $\textbf{z}_{t+T_p}$ in the whole feasible region of control ($T_p$ is the prediction horizon in MPC). (3) \textit{Control compatibility}. The formulation of the WF ROM must be readily incorporated into optimal controller design. 

Given that each grid point has local states denoting the current velocity at that point, the meshed full-state snapshot $\textbf{x}_{t}$ that describes the full WF flow field at time instant $t$ exhibits multiscale spatial physics \cite{kutz2017deep}. Such data provide an opportunity for DNN to make an impact in the modeling and analysis of WF flow fields. For data that have spatial dependencies, convolutional neural networks (CNNs) have been shown to be effective in learning informative high-level features \cite{gupta2014learning}. Another motivation for adopting the DNN for the WF ROM is its flexible architecture that allows for accurately fitting temporal dynamics in the desired form without resorting to hand-designed features \cite{lusch2018deep}. In the following, we will introduce the proposed DNN-based WF ROM in detail.

\subsection{Deep Autoencoder for Reduced-Order WF State Equation} \label{netprob}
To learn an appropriate WF ROM meeting the three requirements, in this paper, we propose a specialized DNN architecture. The proposed DNN belongs to the class of autoencoders. The CNN is utilized as the encoder and decoder layers for its power of handling spatial features to extract the intrinsic coordinates that dominate the WF flow data. A linear dynamic system in the reduced-order space is embedded into the neural network as dynamic bottleneck layers to support the optimal control implementation.

\begin{figure}[tb]
	\centering
		\includegraphics[width=3.7 in]{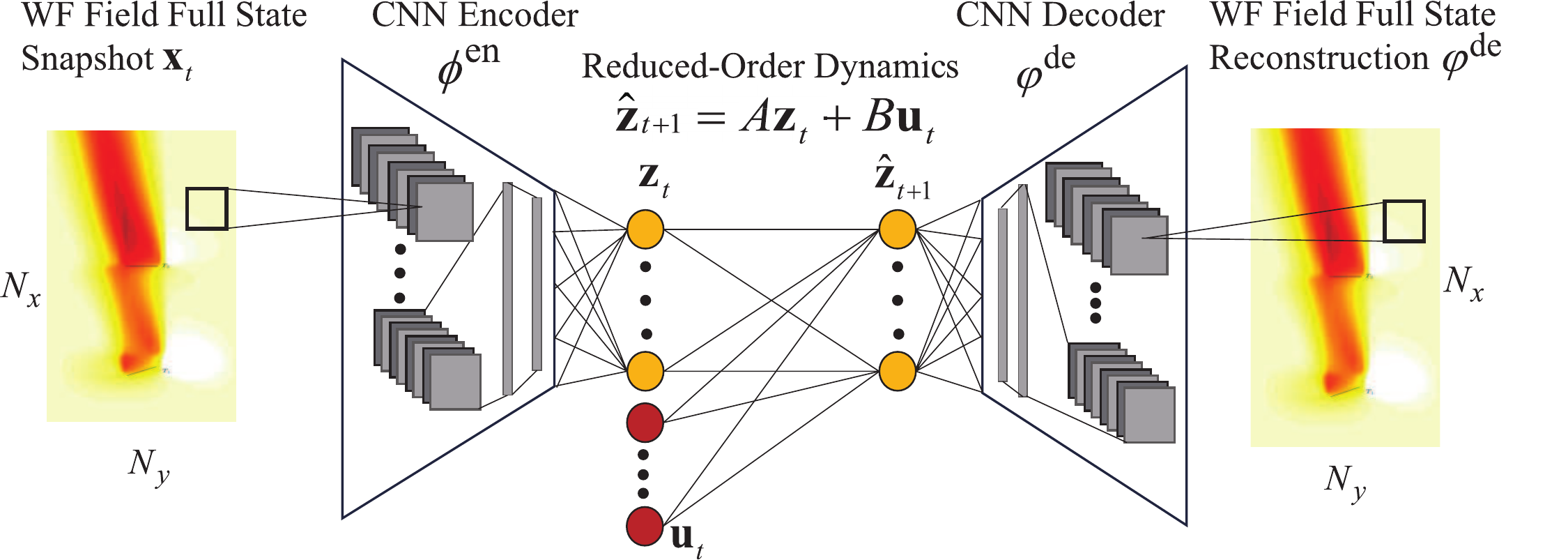}
	\caption{Illustration of the proposed dynamic autoencoder.}
	\label{fig:autoencoder}
\end{figure}

To identify the state equation of WF ROM, a sequence of time snapshots $[\textbf{x}_1,\ldots, \textbf{x}_T]$ and corresponding control inputs $[\textbf{u}_1,\ldots, \textbf{u}_T]$ are required as the training data set; $T$ is the size of the sequential training data set. Fig. \ref{fig:autoencoder} shows a sketch of the proposed dynamic autoencoder architecture.
The high-dimensional wind flow states $\textbf{x}_t$ are compressed through the encoder network $\phi^{\textrm{en}}(\textbf{x}_t)$, which serves as the compression mapping and produces the low-dimensional state $\textbf{z}_t$. We denote the encoder as $\textbf{z}_t = \phi^{\textrm{en}}(\textbf{x}_t)$.
To learn corresponding dynamic in the latent space, the latent state $\textbf{z}_t$ is fed into the designed reduced-order state equation, $\widehat{\textbf{z}}_{t+1} = f^{\textrm{ROM}}(\textbf{z}_t, \textbf{u}_t)$. For optimal control compatibility, we enforce $\widehat{\textbf{z}}_{t+1} = f^{\textrm{ROM}}(\textbf{z}_t, \textbf{u}_t)$ as a linear dynamic form based on Koopman theory \cite{koopman1931hamiltonian}: 
\begin{align}
 \widehat{\textbf{z}}_{t+1} = A\textbf{z}_t+B\textbf{u}_t \label{romdy}
\end{align}
where $A\in \mathbb{R}^{n_z \times n_z}$ $B\in \mathbb{R}^{n_z \times 2N}$ are two global linearization matrices whose parameters are optimized during training along with the neural network parameters. Note that while the dynamics operate on the reduced-order space, it takes untransformed controls into account. The hat symbol $\wedge$ of $\widehat{\textbf{z}}_{t+1}$ denotes the predicted latent states derived from $f^{\textrm{ROM}}$, distinguished from $\textbf{z}_{t+1}$ derived directly from encoder $\phi^{\textrm{en}}(\textbf{x}_{t+1})$.
Reconstruction of the full state from $\textbf{z}_t$ is performed by passing $\textbf{z}_t$ through the decoder network $\varphi^{\textrm{de}}(\textbf{z}_t)$.  

The purpose of the dynamic autoencoder is to map the high-dimensional nonlinear system (\ref{CFDstate}) to a low-dimensional linear system (\ref{romdy}). Koopman theory proves that a nonlinear system can be mapped to an infinite-dimensional linear system through the transformation of system states \cite{koopman1931hamiltonian}, and that the infinite-dimensional linear system can be reduced to a finite-dimensional system if the Koopman invariant subspace is found \cite{koopman1931hamiltonian}. This gives the rationality of the system mapping in this paper. Besides, previous researches have observed the fact that the dominant dynamics of spatial-temporal flow data always evolve in a low-rank subspaces spanned by a few spatial coordinates \cite{kutz2017deep, annoni2017method}. So a linear low dimension system is potential to capture the dominant dynamics of the high dimensional CFD model. However, there are no analytical methods to directly derive such model reduction \cite{kutz2017deep}. The challenge of identifying and representing Koopman invariant subspace provides motivation for the use of deep learning methods for the powerful learning ability. \cite{lusch2018deep} and \cite{morton2018deep} have shown previous success of finding Koopman invariant subspace through deep learning. In the proposed architecture, we enforce the encoder and decoder to extract reduced-order states that evolve linearly in time constrained by (8), attempting to obtain a finite-dimensional approximation of the Koopman operator. The encoder and decoder are the approximations of observable functions that span a Koopman invariant subspace.


The complete loss functions used in training the dynamic autoencoder can be divided into three parts corresponding to the three high-level requirements for the network.
The reconstruction accuracy of the autoencoder is achieved using the following loss:
\begin{align}
  \mathcal{L}_{\textrm{recon}}=\|\textbf{x}_t-\varphi^{\textrm{de}}(\phi^{\textrm{en}}(\textbf{x}_{t}))\|^2_{F}
\end{align}
where $\Vert \cdot \Vert_F$ represents the Frobenius norm, i.e. {\color{black} square root} of the sum of squares of all entries in a tensor, similar to the $2$-norm of a vector. The main purpose of $\mathcal{L}_{\textrm{recon}}$ is to extract the latent states from the high-dimensional data through the encoder and decoder architecture. Thus states at the same time step are utilized to formulate the reconstruction error.


The linear dynamics in WF ROM are constrained through the following loss functions. 
We enforce that $\widehat{\textbf{z}}_{t+m}$ predicted by the reduced-order state equation (\ref{romdy}) to be consistent with $\textbf{z}_{t+m}$ derived from encoding $\phi^{\textrm{en}}(\textbf{x}_{t+m})$:
\begin{align}
  \mathcal{L}_{\textrm{lin}}=\sum\nolimits_{m=1}^{S_\mathrm{p}} \| \widehat{\textbf{z}}_{t+m} - \textbf{z}_{t+m} \|^2_2
\end{align}
where $S_{\mathrm{p}}$ is a hyperparameter for how many steps to check in the loss function. $\widehat{\textbf{z}}_{t+m}$ is derived from evolving  $\textbf{z}_{t}$ through the dynamic equation (\ref{romdy}) $m$ steps forward in time.

To enable optimal control in ROM, accurate future state prediction is necessary. To this end, the state equation must capture accurate dynamic relations between control inputs and reduced-order states. This requirement is reflected through the following loss:
\begin{align}
  \mathcal{L}_{\textrm{pre}}= \sum\nolimits_{m=1}^{S_{\mathrm{p}}} \|\textbf{x}_{t+m}-\varphi^{\textrm{de}}(\widehat{\textbf{z}}_{t+m}) \|^2_F 
\end{align}

\begin{figure}[h]
	\centering
		\includegraphics[width=3.7 in]{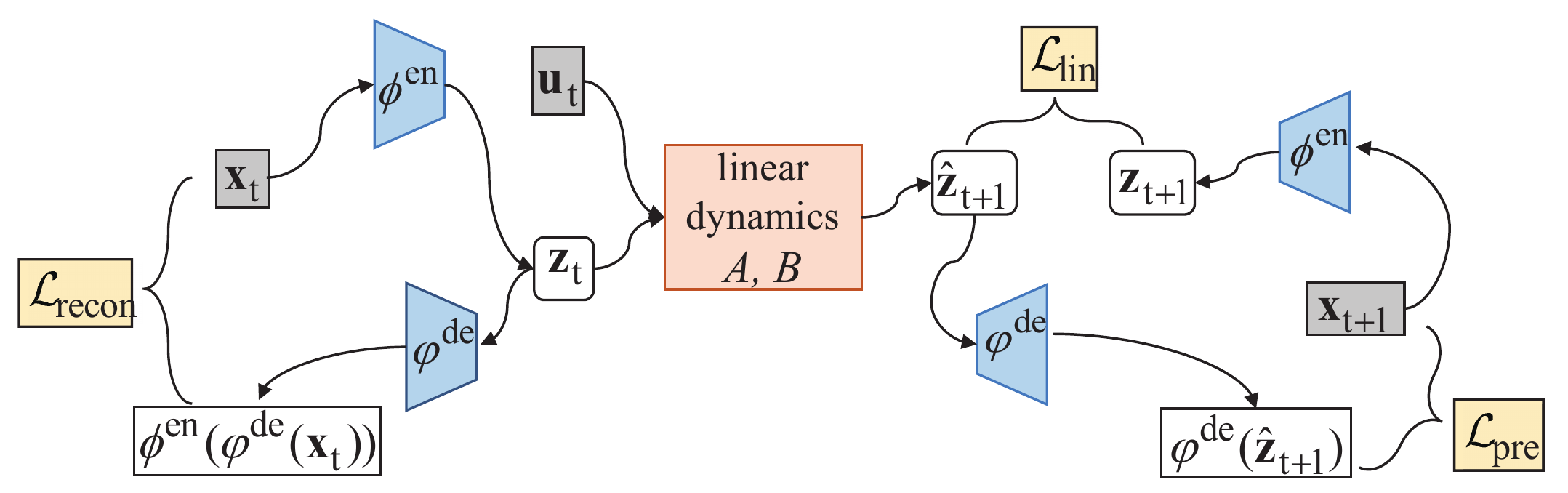}
	\caption{ Data flow in the training procedure of the dynamic autoencoder. In practice, $\widehat{\textbf{z}}$ is evolved through (\ref{romdy}) repeatedly to calculate $\mathcal{L}_{\textrm{pre}}$ and $\mathcal{L}_{\textrm{lin}}$.}
	\label{fig:aeflowchart}
\end{figure}

The overall loss function to train the dynamic autoencoder is the sum of the above three parts:
\begin{align}
  \mathcal{L}=\mathcal{L}_{\textrm{recon}}+ \beta \mathcal{L}_{\textrm{pre}}+ \alpha \mathcal{L}_{\textrm{lin}}
\end{align}
We apply $L_2$ regularization to the network parameters. The weights $\alpha$ and $\beta$ are hyperparameters to combine the three losses. The overall flowchart of the training is illustrated in Fig. \ref{fig:aeflowchart}.
The encoder consists of the widely used CNN structure, ResNet convolutional layers \cite{he2016deep}, while the decoder inverts all operations performed by the encoder. The detailed network parameters are elaborated in Section \ref{caseset}.

A sketch of the data flow in the training procedure is shown in Fig. \ref{fig:aeflowchart}. In practice, $\textbf{z}_t$ is evolved through the linear dynamic equation (\ref{romdy}) $S_p$ steps forward in time to calculate $\mathcal{L}_{\textrm{pre}}$ and $\mathcal{L}_{\textrm{lin}}$. To learn a reduced order system for optimal control application, the three parts of the loss function are integral and necessary. Minimizing $\mathcal{L}_{\textrm{recon}}$ enforces the mapping is invertible so that the dimension reduction is valid, which is the typical loss function for a standard autoencoder. Minimizing $\mathcal{L}_{\textrm{pre}}$ enforces that the derived dynamic system can accurately simulate the time evolution of original CFD model. Minimizing $\mathcal{L}_{\textrm{lin}}$ enforces the latent state uniqueness. Without enforcing the consistency between $\widehat{\textbf{z}}_{t+m}$ and $\textbf{z}_{t+m}$, (\ref{romdy}) may produce a latent state from which reconstruction of $\textbf{x}$ is possible but that is not a valid encoding (i.e. the corresponding latent state derived by the encoder).

In the proposed dynamic autoencoder, the parameters of the matrix A and B are trained with the encoder and decoder parameters by minimizing the loss function (12). The state reconstruction and system dynamics are simultaneously considered in the training procedure. Because the latent states for purely reconstruction can not cover important system dynamic information \cite{kutz2017deep} \cite{ lusch2018deep}, such integration is necessary for the model reduction of dynamic systems.

%

\subsection{Reduced-Order WF Power Output Equation} \label{netprob}
To control the WF dynamic power production based on WF ROM, the power output equation in the reduced-order space, $\textbf{P}_t = f^{\textrm{pow}}(\textbf{z}_t, \textbf{u}_t) \label{euq:pow}$, is required to be identified in the form that can be incorporated into optimal control. The reduced-order states $[\textbf{z}_1,\ldots, \textbf{z}_T]$ derived from the trained encoder, control inputs $[\textbf{u}_1,\ldots, \textbf{u}_T]$ and corresponding wind farm power productions $[\textbf{P}_1,\ldots, \textbf{P}_T]$ are used as the training data set.

Note that although this output relation is straightforward because no dynamics exist, it needs to fulfill the control compatibility requirement.
Directly applying the nonlinear output equation in control implementation will cause great difficulty for optimization \cite{watter2015embed}. We circumvent this problem by imposing a locally linear form on this output equation in the training process:
\begin{align}
  \textbf{P}_t = C_t \textbf{z}_t+D_t \textbf{u}_t+o_t \label{rompower} 
\end{align}
where $C_t\in \mathbb{R}^{N \times n_z}$, $D_t\in \mathbb{R}^{N \times 2N}$ are local Jacobians, and $o_t\in \mathbb{R}^{N}$ is the offset. These matrices $C_t , D_t$ and $o_t$ formulate the local linearization of the nonlinear output equation and are determined at each time step as functions of the current reduced-order states and control inputs. 

\begin{figure}[h]
	\centering
		\includegraphics[width=2.6 in]{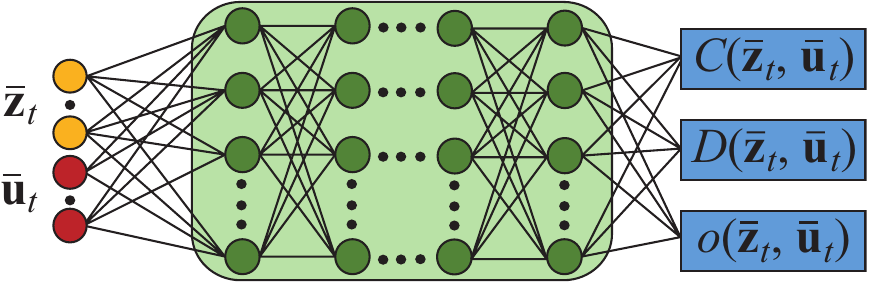}
	\caption{Schematic of $\Psi(\cdot)$ to train the local linearization matrices.}
	\label{fig:rewardNN}
\end{figure}

Therefore, instead of directly learning the output equation, a separate neural network is designed to learn the linearization matrices, as shown in Fig. \ref{fig:rewardNN}.
This neural network is a fully-connected one and maps from the current reduced-order states and control inputs $(\overline{\textbf{z}}_t, \overline{\textbf{u}}_t)$ to the corresponding local linearization matrices.
We denote this neural network as $ [C_t, D_t, o_t]  = \Psi(\overline{\textbf{z}}_t,\overline{\textbf{u}}_t)$. To ensure the accuracy of the local linearization, in the training process, both $(\overline{\textbf{z}}_t, \overline{\textbf{u}}_t)$ and its
neighborhood point $(\textbf{z}_t , \textbf{u}_t)$ are required. To this end, we adds a random noise of small amplitude to the original $(\textbf{z}_t , \textbf{u}_t)$ in the training data set to generate $(\overline{\textbf{z}}_t, \overline{\textbf{u}}_t)$. $(\overline{\textbf{z}}_t, \overline{\textbf{u}}_t) = (\textbf{z}_t , \textbf{u}_t)+ \mathcal{N}(0, \epsilon)$.
 Then, the loss function for training $\Psi$ is as follows:
\begin{align}
  \mathcal{L}_{\mathrm{pow}}=\|  \textbf{P}_t - C(\overline{\textbf{z}}_t, \overline{\textbf{u}}_t)\textbf{z}_t+D(\overline{\textbf{z}}_t, \overline{\textbf{u}}_t)\textbf{u}_t+o(\overline{\textbf{z}}_t, \overline{\textbf{u}}_t) \|^2_2  \label{rompowertrain}
\end{align}

After training, the local linearization matrices can be directly derived through $\Psi$. In this way, in the optimal control algorithm, the locally linear form (\ref{rompower}) can be readily formulated based on the derived linearization matrices.

Through training the proposed two networks, the reduced-order state equation, Eq. (\ref{romdy}), and the power output equation, Eq. (\ref{rompower}), can be finally derived and compose the WF ROM.

\section{The Control Framework Formulation} \label{control}
When the training is done, the obtained WF ROM can be used for online control applications. In this section, a novel control algorithm is introduced so that the WF can optimally control its dynamic power production to track the power reference given by the power system. The control framework explicitly incorporates the WF ROM to solve the optimal control problem efficiently. 

\begin{figure}[h]
	\centering
		\includegraphics[width=3.7 in]{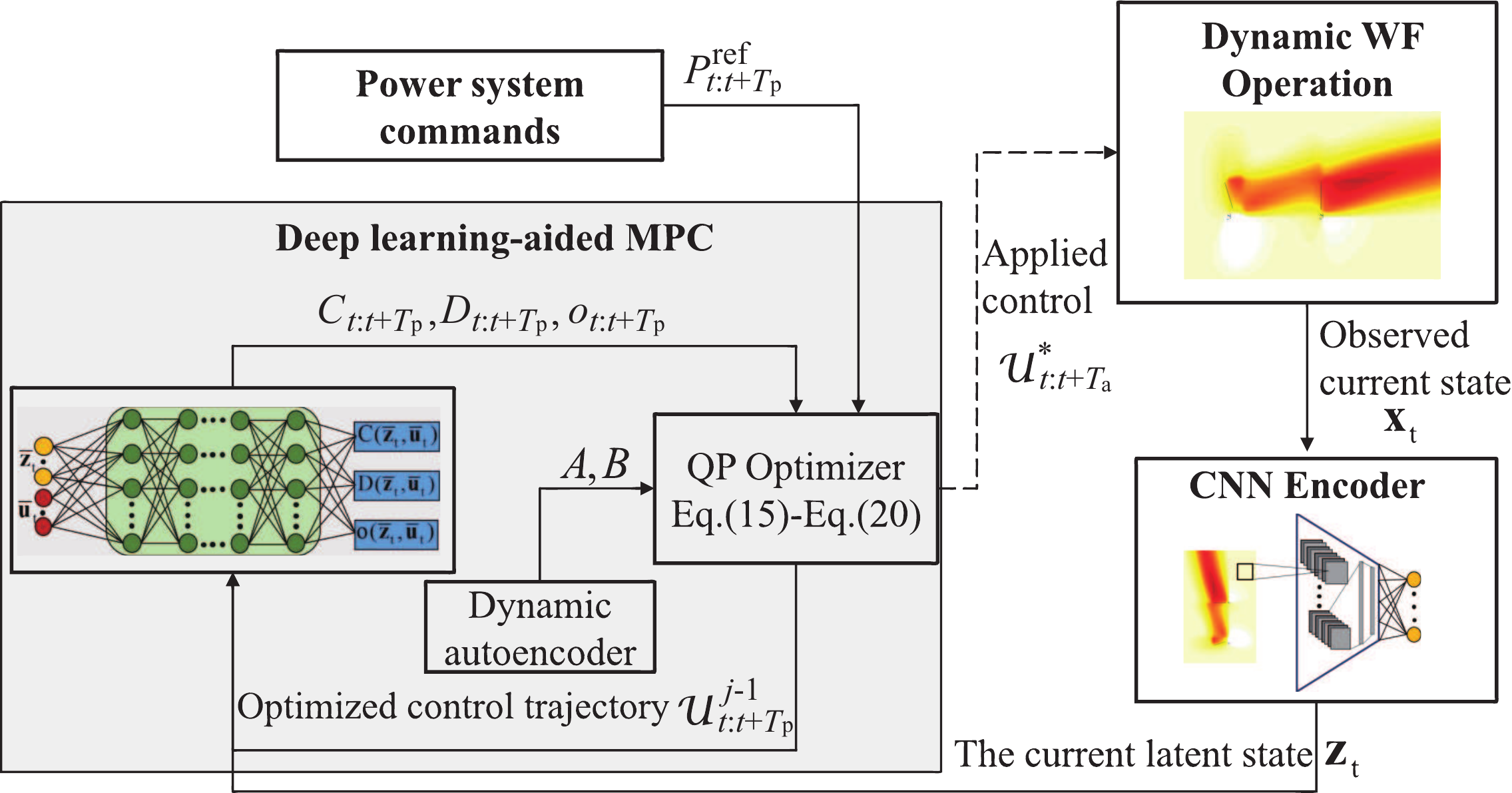}
	\caption{Overall scheme of the proposed deep learning-aided WF MPC.}
	\label{fig:wake}
\end{figure}
The overall WF control framework is shown in Fig. \ref{fig:wake}. At each control period, the current full-state WF snapshot $\textbf{x}_t$ can be observed from CFD simulators or the wind farm field radar observer \cite{doekemeijer2016enhanced, doekemeijer2017ensemble}. Then, the current reduced-order state $\textbf{z}_t$ is derived by passing $\textbf{x}_t$ through the trained encoder $\phi^{\textrm{en}}(\textbf{x}_{t})$. The reduced-order dynamic equation (\ref{romdy}) is globally linear. The power equation is nonlinear, but its locally linear formulation has been learned by (\ref{rompower}). Based on the current state $\textbf{z}_t$, the optimal control problem on the learned WF ROM gives us a locally linear quadratic formulation. 
After the optimal control trajectory is found, the first $T_a$ control steps are applied to the WF to follow the power reference. When new system states are observed after $T_a$ time, a new control sequence can be found by repeating this procedure using the last derived control trajectory as initialization.

After the completion of training in Section \ref{network}, the optimal control is conducted with the derived WF ROM independently from CFD simulation. The optimization is performed entirely in the reduced-order space. Thus, the control complexity and the computational burden is greatly reduced compared to directly controlling the full-order system. 
The detailed control algorithm is elaborated in Algorithm 1.

At each iteration $j$ in Algorithm 1, the quadratic optimization problem  $\mathcal{P}_{t}^{j}$ is formulated as follows:
\begin{align}
 \min\limits_{{\mathcal{U}}_{t:t+T_\mathrm{p}}^j} \sum_{k=t}^{t+T_\mathrm{p}} [(P^{\mathrm{ref}}_{k}-P^{\textrm{WF}}_{k})^\textrm{T} Q (P^{\mathrm{ref}}_{k} - P^{\textrm{WF}}_{k}) + (\textbf{u}_{k}^j-\textbf{u}_{k-1}^j)^\textrm{T} R  (\textbf{u}_{k}^j-\textbf{u}_{k-1}^j)],    \label{obj}
 \end{align}
\vspace{-30pt}
\begin{align}
\text{s.t.:} \hspace{15pt} \textbf{z}_{k+1}  =   A\textbf{z}_{k}+B\textbf{u}_{k}^j, \hspace{30pt} k=t,\ldots,t+T_{\mathrm{p}},  \label{opt:dyn}  \hspace{50pt} \\
 \textbf{P}_{k}  =  C(  \textbf{z}_k^{j-1},  \textbf{u}_k^{j-1} )\textbf{z}_{k} +D(\textbf{z}_k^{j-1}, \textbf{u}_k^{j-1}) \textbf{u}_{k}^j + o(\textbf{z}_k^{j-1}, \textbf{u}_k^{j-1} ), \hspace{10pt} k =  t,\ldots,t+T_{\mathrm{p}},   \label{opt:pow}\\
\textbf{u}^{\textrm{min}}  \leq \textbf{u}_{k}^j\leq \textbf{u}^{\textrm{max}}, \hspace{30pt} k=t,\ldots,t+T_{\mathrm{p}},  \label{opt:bound} \hspace{50pt} \\
-\triangle\textbf{u}^{\textrm{max}}  \leq \textbf{u}_{k}^j -\textbf{u}_{k-1}^j\leq \triangle\textbf{u}^{\textrm{max}}, \hspace{15pt} k=t,\ldots,t+T_{\mathrm{p}},  \hspace{35pt} \label{opt:delt} \\
\textbf{u}_k^{j-1} - \varepsilon  \leq \textbf{u}_{k}^j\leq \textbf{u}_k^{j-1} + \varepsilon, \hspace{15pt} k=t,\ldots,t+T_{\mathrm{p}}, \hspace{45pt} \label{opt:simil}
\end{align}
where the subscript $k$ represents the time step, and $t$ represents the current time step. $j$ is the iteration number in Algorithm 1.
$T_p$ denotes the prediction horizon. The decision variable $\mathcal{U}_{t:t+T_p}^j$ is an aggregated vector of the time series of control inputs over the prediction horizon:
\begin{align}
\mathcal{U}_{t:t+T_p}^j \overset{\textrm{def}}{=} [\textbf{u}_t^j, \hspace{5pt} \textbf{u}_{t+1}^j, \hspace{5pt} \ldots, \hspace{5pt} \textbf{u}_{t+T_p}^j]
\end{align}

The objective function (\ref{obj}) tries to find a distribution of control signals among the turbines such that the tracking error and dynamical input variation are minimized. The first part minimizes the deviation between the power reference $P^{\textrm{ref}}$ and the WF power production $P^{\textrm{WF}}$ to enable the WF to provide AGC service. $P^{\textrm{WF}}_k = \sum_{i=1}^{N}P^i_k$. The second part penalizes control actuator costs to prevent oscillations in the input commands. $Q$ and $R$ are corresponding weighting matrices.

The constraints (\ref{opt:dyn}) and (\ref{opt:pow}) derived from Section \ref{network} describe the WF ROM. The dynamic equation (\ref{opt:dyn}) is linear. The power equation (\ref{opt:pow}) is also linear at each time step once $[C_k, D_k, o_k]$ are determined. Even though  varying over the prediction horizon, they are essentially constant matrices unrelated to the decision variable. These local linearization matrices are determined by the last state and control trajectory of $(\textbf{z}_k^{j-1}, \textbf{u}_k^{j-1} )$  through the second trained DNN $\Psi (\cdot)$ in Fig. \ref{fig:rewardNN}. The last state and control trajectory are derived by evolving the initial latent state $\textbf{z}_t$ through the learned dynamic equation (\ref{romdy}) based on the last optimization control trajectory $\mathcal{U}_{t:t+T}^{j-1}$ calculated by solving $\mathcal{P}_{t}^{j-1}$ at last iteration. 

The constraints (\ref{opt:bound}) and (\ref{opt:delt}) ensure a feasible region of a practical WF control signal, including the upper and lower bounds, and the ramping restrictions. The constraint (\ref{opt:simil}) enforces that the decision control trajectory $[\textbf{u}_t^j \ldots \textbf{u}_{t+T_p}^j]$ is close to the last optimization results $[\textbf{u}_t^{j-1} \ldots \textbf{u}_{t+T_p}^{j-1}] $ to ensure the accuracy of the locally linear power equation (\ref{opt:pow}).
$\varepsilon$ is a small number. Without this constraint, the calculated control trajectory may be far from the one used to determine the local linearization matrices in (\ref{opt:pow}), and its accuracy may be affected. The increased error of constraint (\ref{opt:pow}) deteriorates the convergence of Algorithm 1.

A computational complexity analysis of the optimization is also provided. The main computation effort is to solve the optimization problem. If directly controlling the CFD model based on adjoint method, in each iteration, the time complexity for evaluating the derivatives of the objective function to control variables is $O(n_{x}^{3})$ by CFD simulation \cite{munters2018dynamic, vali2019adjoint}, where $n_{x}$ is the number of CFD model states. For the proposed method, in each iteration of Algorithm 1, the complexity for evaluating the derivatives is approximately $O(n_{z}^{2} + n_{z}n_{u})$, where $n_{u}$ is the number of control inputs and $n_{z}$ is the number of ROM states. Because the linear system and the local linearization matrixes in (\ref{opt:dyn}) and (\ref{opt:pow}) are already obtained from the trained networks, no matrix inverse is required to evaluate the derivatives and the corresponding complexity is thus square. Because $n_{z}, n_{u}\ll n_{x}$, the complexity is reduced by the formulation of WF ROM. 

\begin{algorithm}[tb]
\caption{Deep Learning-Aided MPC for WF control}
{\bf Input:} 
Current full-order WF state $\textbf{x}_{t}$. Last optimal control trajectory $\mathcal{U}_{(t-1):(t-1)+T_{\mathrm{p}}}^{\ast}$. Error tolerance $\zeta$ \\
{\bf Output:}  The optimal control trajectory $\mathcal{U}_{t:t+T_{\mathrm{p}}}^{\ast}$
\begin{algorithmic}[1]
\State \textbf{Initialization}: Pass $\textbf{x}_{t}$ through the encoder $\phi^{\textrm{en}}$ to obtain $\textbf{z}_{t}$. Initialize the control trajectory $\mathcal{U}_{t:t+T_{\mathrm{p}}}^0 = \mathcal{U}_{(t-1):(t-1)+T_{\mathrm{p}}}^{\ast}$. Initialize $j = 0$.
 \Repeat
    \State $j \leftarrow j+1$
    \For{$k=t$ to $t+T_p$}
        \State Obtain the coefficient matrices $C_k, D_k, o_k$ from the trained network $\Psi$ with $\textbf{z}_{k}$ and $\textbf{u}_{k}^{j-1}$.
        \State Evolving $\textbf{z}_{k}$ forward in time through the reduced-order dynamic equation (9) with $\textbf{u}_{k}^{j-1}$.
    \EndFor
    \State \textbf{end for}
    \State Solve the quadratic programming problem $\mathcal{P}_{t}^{j}$ (\ref{obj})-(\ref{opt:simil}), obtain optimized control trajectory $\mathcal{U}_{t:t+T}^{j}$
\Until{($\|\mathcal{U}_{t:t+T_{\mathrm{p}}}^j - \mathcal{U}_{t:t+T_{\mathrm{p}}}^{j-1} \| < \zeta$)}  \\
\Return the optimal control trajectory $\mathcal{U}_{t:t+T_{\mathrm{p}}}^{\ast}$
\State Apply the first $T_a$ control actions $\mathcal{U}_{t:t+T_{\mathrm{a}}}^{\ast}$ to the real WF
\end{algorithmic}
\end{algorithm}


\section{Case Study}\label{case}
In this section, we verify the performance of the proposed deep learning-aided MPC. The WF simulation is based on the NS equation-based open-source simulator WFSim developed by the Delft University of Technology \cite{boersma2018control}. WFSim has been widely employed in WF studies \cite{boersma2017tutorial, vali2018model, zhao2020cooperative}. Its accuracy have been validated by previous studies \cite{boersma2018control} \cite{vali2019adjoint}. For simulator validation, readers are referred to previous studies \cite{boersma2018control} and \cite{vali2019adjoint} which have demonstrated that WFSim is able to capture the dominant dynamics of the WF wakes. We also further validate its accuracy based on realistic measurement data of Garrad Hassan, which can be found online\footnote{The data and validation description has been uploaded to https://github.com/kkxchen/WF-deep-learning-aided-MPC-supplementary-material for easy access}. The DNN is constructed in the \emph{Python} environment utilizing the \emph{TensorFlow} framework.

\subsection{Case Settings} \label{caseset}
A layout of a $3 \times 3$ WF is considered and simulated with WFSim. The typical NREL 5 MW Type III WT is assumed to be installed.

The rotor diameter D $=126$ m \cite{jonkman2009definition}. WTs have 5D spacing in the streamwise direction and 3D in the spanwise direction. The simulation parameters are set as follows \cite{boersma2018control, munters2018dynamic}.
The air density $\rho=1.2\textrm{kg}/\textrm{m}^{3}$. The upper and lower bound of the thrust coefficient are set as $0.1 \leq C_{T}^\prime \leq 2 $. The value for $C_{T}^{\prime \textrm{max}} =2$  corresponds to the Betz-optimal value. $C_{T}^{\prime \textrm{min}} = 0.1$ indicating that we do not allow turbines to shut down completely. The turbine electrical dynamics whose time scale is in millisecond are not considered in this study, which is common practice in WF level researches \cite{boersma2017tutorial}. The ramping limit $\triangle C_{T}^{\prime \textrm{max}}$ is $0.2/\textrm{s}$ such that turbines can not de- and uprate instantaneously. The control limits for yaw angles are as follows. $ -25^\textrm{o}\leq \gamma \leq 25^\textrm{o}$. $\triangle \gamma^{\textrm{max}} = 0.3^\textrm{o}/\textrm{s}$. In the simulation, the simulated wind field is $2520\times1560 \textrm{m}^2$, with a grid of $100\times55$ cells ($N_x \times N_y$). The corresponding cell size of the discretization is $25.2m \times 28m$. The sample period is 1 s. The incoming wind speed in the experiment is $10\textrm{m}/\textrm{s}$. The wake-loss-heavy situation is studied when the wind
direction aligns with the row of turbines.

Analogous to frequency sweeps in system identification \cite{morton2018deep}, we excite the system with sinusoidal control inputs with random frequencies during simulation to generate data. Every 500 s of simulation, the control input frequencies of all turbines are altered and randomly selected. 
Time snapshots of the current $\textbf{x}_t$, $\textbf{u}_t$ and $\textbf{P}_{t}$ are stored every one second of simulation. In total, the training set contains $10^{5}$ sequential snapshots of the dynamic WF simulations. To test the generalization ability of the trained models, the validation data set contains 5000 samples generated in the same way. 

\begin{table}[th]
\centering
\caption{Network Structure of the Encoder}
\label{Tab1}
\begin{tabular}{c c c }
\hline
  Block name & Output size & Layer parameters \\
\hline
  Res Block.1 & $50 \times 28$ &\fontsize{9pt}{10pt} $\left[\begin{smallmatrix}  1 \times 1, 64 \\  3 \times 3, 64 \\  1 \times 1, 256 \end{smallmatrix}\right] \times 3$ \\
  Res Block.2 & $25 \times 14$ &\fontsize{9pt}{10pt}  $\left[\begin{smallmatrix} 1 \times 1, 128 \\  3 \times 3, 128 \\  1 \times 1, 512 \end{smallmatrix}\right] \times 3$ \\
  Res Block.3 & $13 \times 7$ &\fontsize{9pt}{10pt}  $\left[\begin{smallmatrix} 1 \times 1, 256 \\  3 \times 3, 256 \\  1 \times 1, 1024 \end{smallmatrix}\right] \times 3$ \\
  Dense & $n_z=20$ & $\begin{smallmatrix} 3\hspace{2pt}\textrm{fully} \hspace{2pt}\textrm{ connected} \hspace{2pt} \textrm{layers} \\ 2000, 1200, 500 \end{smallmatrix}$  \\
\hline
\end{tabular}
\end{table}

For the dynamic autoencoder, the WF full-state dimension is $100\times55\times2$, i.e. 11000. The reduced-order space dimension $n_z$ is set to 20. The architecture of the encoder is shown in Table \ref{Tab1}. The widely used CNN structure, residual convolutional blocks \cite{he2016deep}, with ReLU activations construct the encoder. The decoder inverts all operations performed by the encoder. For the loss function, $S_{p} = 50$. $\alpha = 300$ and $\beta = 1/S_{p}$ are set to balance the three losses so that they are at about the same order of magnitude.

\subsection{Test of the Accuracy of the Learned WF ROM} \label{accuracy}
We compare the prediction accuracy of the proposed DNN method with the mainstream classic model reduction method, dynamic mode decomposition with control (DMDc) method. The DMDc method has been recently used to construct the WF ROM by \cite{annoni2017method, annoni2016wind}. DMDc tries to devise a low-rank linear system that best approximates the system dynamics. DMD is based on singular value decomposition (SVD) to derive the interpretable modes that characterize spatiotemporal flow data. The dimension reduction is done by a linear projection to the modes provided by SVD.
To test the prediction accuracy, WF ROMs derived from the DMDc and DNN approaches are used to predict the future state and power production trajectories. 
Then, the predicted trajectories are tested against the WFSim simulation results.

\begin{figure}[b]
\centering
 \begin{minipage}[t]{0.4\linewidth}
    \centerline{\includegraphics[width=2.4 in]{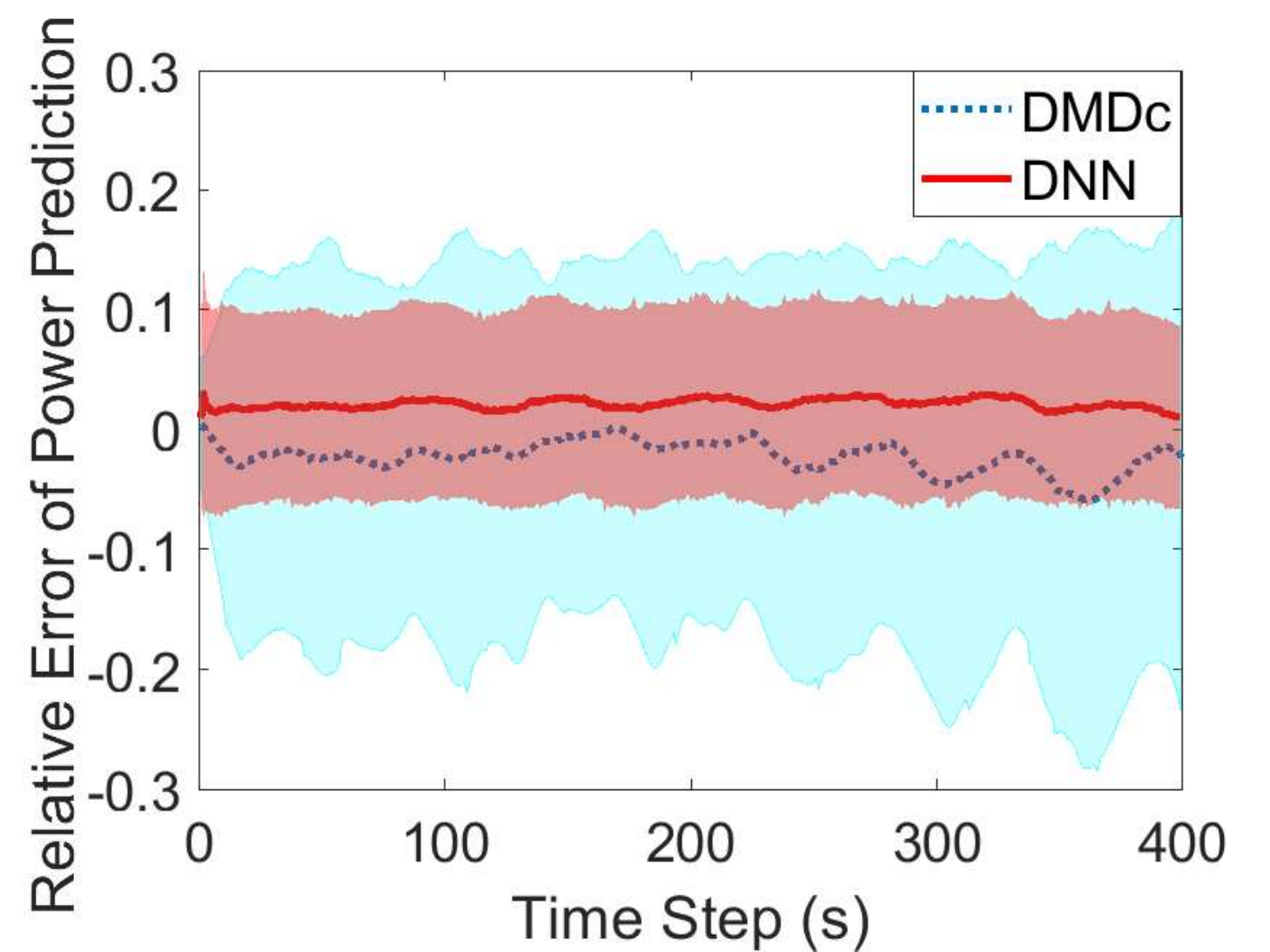}}
     \centerline{\footnotesize{(a) relative error in the training set}}
 \end{minipage}%
 \begin{minipage}[t]{0.4\linewidth}
    \centerline{\includegraphics[width=2.4 in]{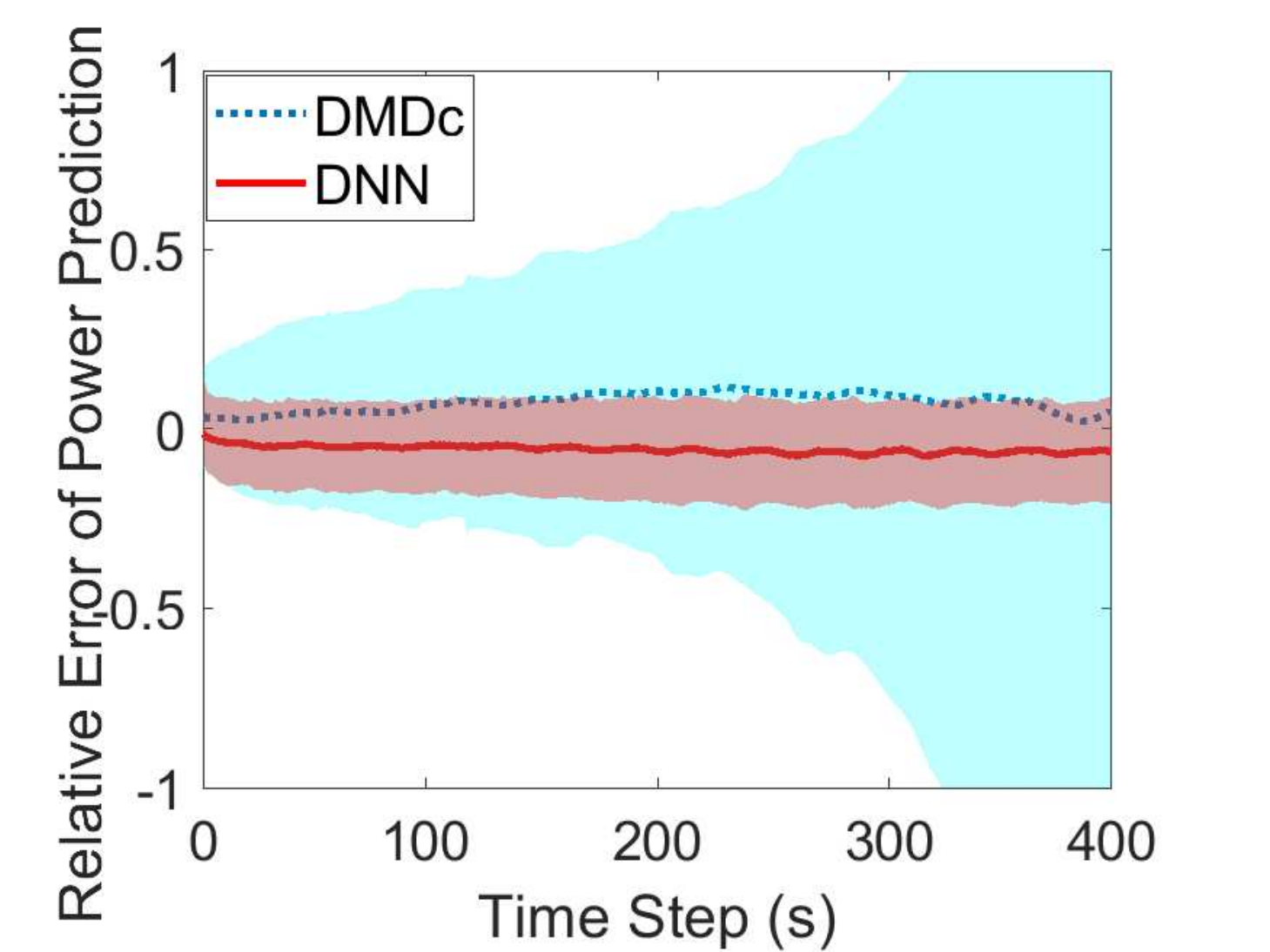}}
     \centerline{\footnotesize{(b) relative error in the test set}}
  \end{minipage}
 \caption{Relative power prediction error. Solid lines represent the mean prediction error across 300 tests, while the shaded regions correspond to one standard deviation about the mean.}
\label{fig:power_error}

 \end{figure}

300 tests are conducted for each method with a power prediction of 400 time steps forward in each test. Fig. \ref{fig:power_error} shows the relative error against the full-order WFSim simulation both in the training data set and testing data set. Solid lines represent the mean prediction error across 300 tests, while the shaded regions correspond to one standard deviation about the mean. Both DMDc and the proposed DNN obtain a good accuracy for reconstructing the already observed time evolution in the training data set. However, the DMDc performance degrades considerably in predicting the test data. This occurs because the linear nature of DMDc limits its generalization ability. 
Besides, SVD encounters difficulty in handling multiscale spatial features \cite{kutz2017deep}. In fact, DMDc can be viewed as a simplified autoencoder with a single unactivated fully connected layer as the encoder. In contrast, the proposed DNN is able to generate stable predictions even in the test data. Its relative error of power prediction remains less than $0.1$. The CNN inherently handles spatial features and the deep autoencoder provides nonlinear coordinations.

\begin{figure} [tb]
	\centering
		\includegraphics[width=3.7 in]{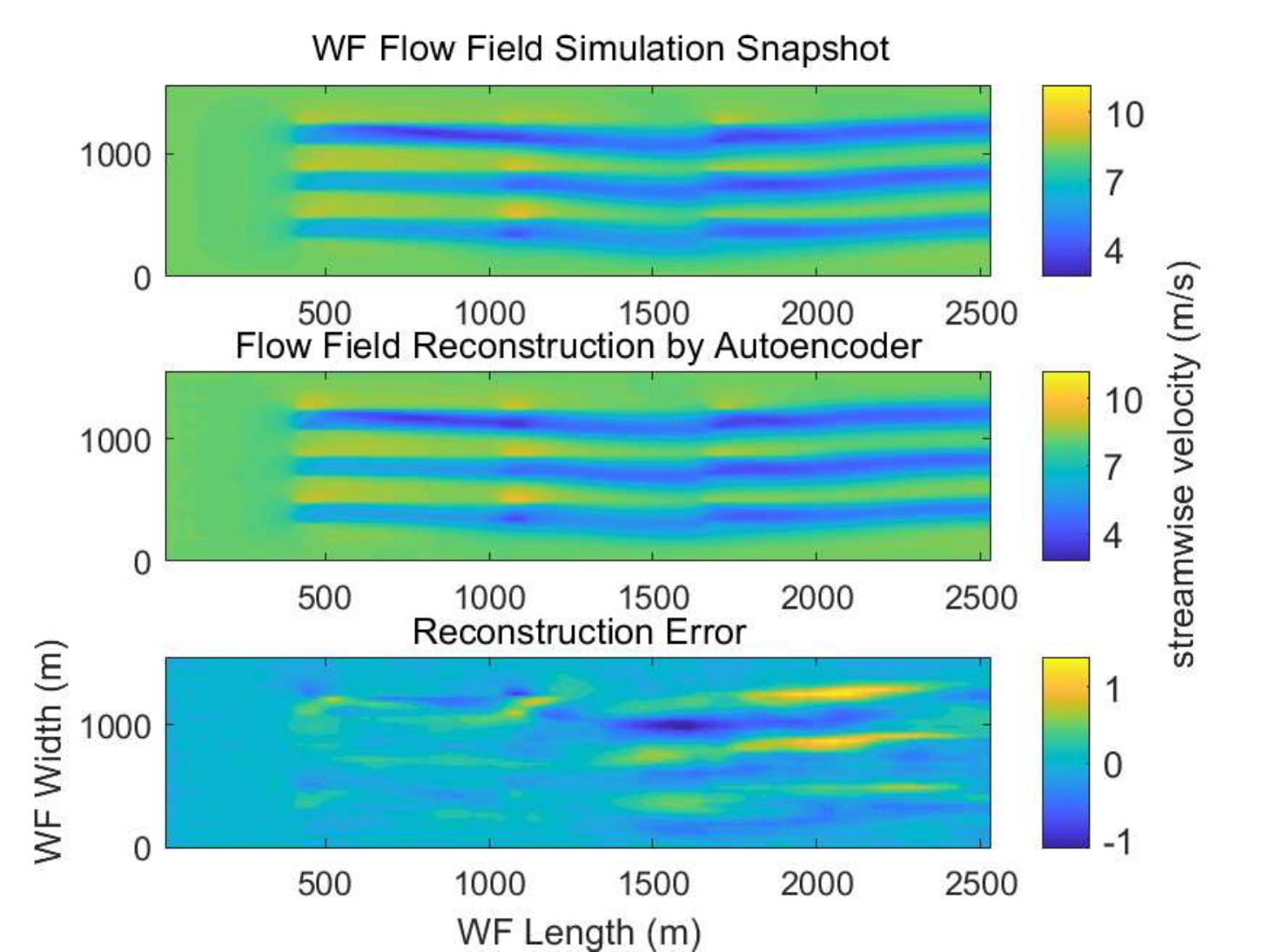}
\vspace{-4pt}	
\caption{Comparison of wind velocity field snapshots by WFSim simulation and dynamic autoencoder prediction.}
\vspace{-12pt}
	\label{fig:WFpic}
\end{figure}

The field visualization of streamwise velocity at a snapshot is shown in Fig. \ref{fig:WFpic}. The top figure is simulated by the NS-based WFSim. The middle one is predicted and reconstructed by the trained dynamic autoencoder. The bottom is the difference between the two images.
In the entire training data set, the mean absolute reconstruction error of the whole $100\times55$ wind velocity field is $0.3\textrm{m}/\textrm{s}$ with the average relative error of $3.26\%$ and maximal relative error of $9.02\%$. In the test data set, the mean reconstruction error is $0.69\textrm{m}/\textrm{s}$ with the average relative error of $6.41\%$. The results show the errors are much less than the simulated ground truth.

\subsection{Deep Learning-Aided MPC Performance} \label{9wtcase}

The widely used proportional distribution (ProD) in current studies is illustrated as the benchmark \cite{sorensen2005wind, van2017active}. At each control step, each turbine estimates its available power capabilities based on the measured wind velocity and sends it to the WF controller. Then, the total power demand is distributed proportionally to the available power capabilities of each turbine \cite{sorensen2005wind}. No wake interaction is considered. 
The greedy power $P_{\textrm{greedy}}$ is defined as the time-averaged WF power production with local greedy settings, where $C_{T}^\prime =2$ and $\gamma = 0^{\textrm{o}}$ for all turbines in the WF \cite{boersma2019constrained}. $P_{\textrm{greedy}}$ is considered as the maximal trackable power by previous WF control schemes that do not consider the wake effect, such as ProD \cite{sorensen2005wind, van2017active}.

\begin{figure}[h]
	\centering
	\includegraphics[width=3.6 in]{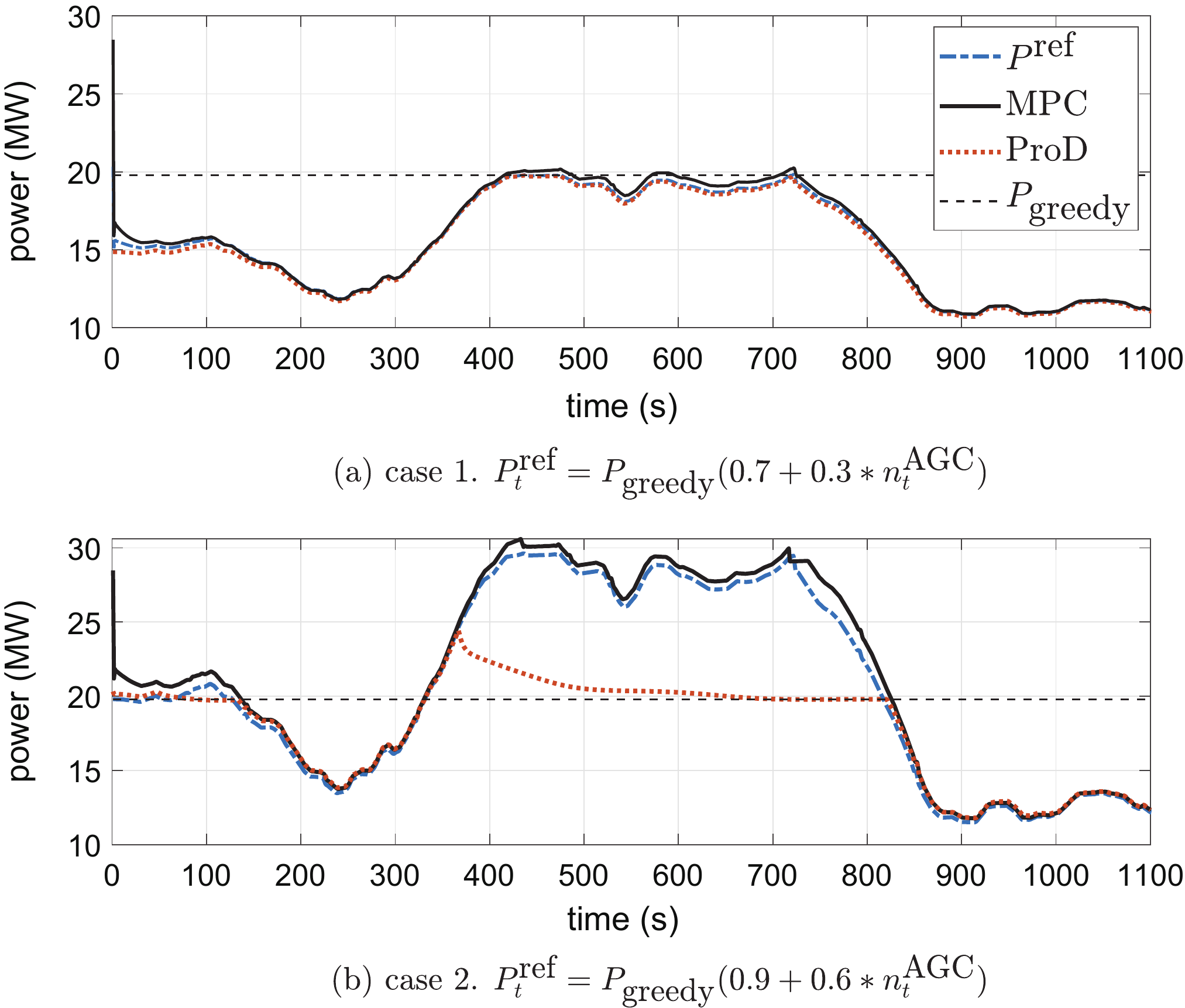}
	\caption{WF power tracking performance.}
	\label{fig:AGCgreedy}
\end{figure}

The commonly used power reference signal for standard AGC qualification, the RegD signal coming from the Pennsylvania, New Jersey, Maryland (PJM) electric market \cite{pilong201712}, is applied to test the control performance, as shown in Fig. \ref{fig:AGCgreedy}.
Two different scenarios are evaluated. For the top figure case, the time-varying power reference is set to $ P^{\textrm{ref}}_t=P_{\textrm{greedy}}(0.7+ 0.3* n_t^{\textrm{AGC}})$. $n_t^{\textrm{AGC}}$ is the RegD test signal that is normalized to $\pm 1$.  Therefore, this power reference will never exceed $P_{\textrm{greedy}}$.
For the bottom figure case, the power reference is $ P^{\textrm{ref}}_t=P_{\textrm{greedy}}(0.9+ 0.6* n_t^{\textrm{AGC}})$.
For a period, more power is demanded than $P_{\textrm{greedy}}$.

As shown in Fig. \ref{fig:AGCgreedy}, if the power reference is lower than $P_{\textrm{greedy}}$, ProD can provide great tracking performance without considering the wake interaction. However, ProD fails when the power reference is higher. Interestingly, as shown in Fig. \ref{fig:AGCgreedy} (b), from $t =320$ s to $t=370$ s, the WF power produces more than $P_{\textrm{greedy}}$ with ProD control, which occurs because wakes of upstream turbines are not yet fully developed. When the wake arrives at downstream turbines, the available wind power decreases and the power production converges to $P_{\textrm{greedy}}$ from $t =370$ s to $t=600$ s. 

In contrast, power tracking can be ensured in both cases by the proposed deep learning-aided MPC. The power regulation capacity is increased compared to ProD. This occurs because the proposed MPC can take into consideration the dynamic wake interactions among WTs in a time horizon $T_p$. If $T_p$ is long enough to cover the wake propagation time within the given WF, MPC can coordinate each turbine operation dynamically to increase the power production capacity while maintaining the power tracking performance. 
For this case, the prediction horizon $T_p =250$s, the receding horizon $T_a =30$s. $\varepsilon=0.02$.

\begin{figure}[h]
	\centering\vspace{-10pt}
		\includegraphics[width=3.0in]{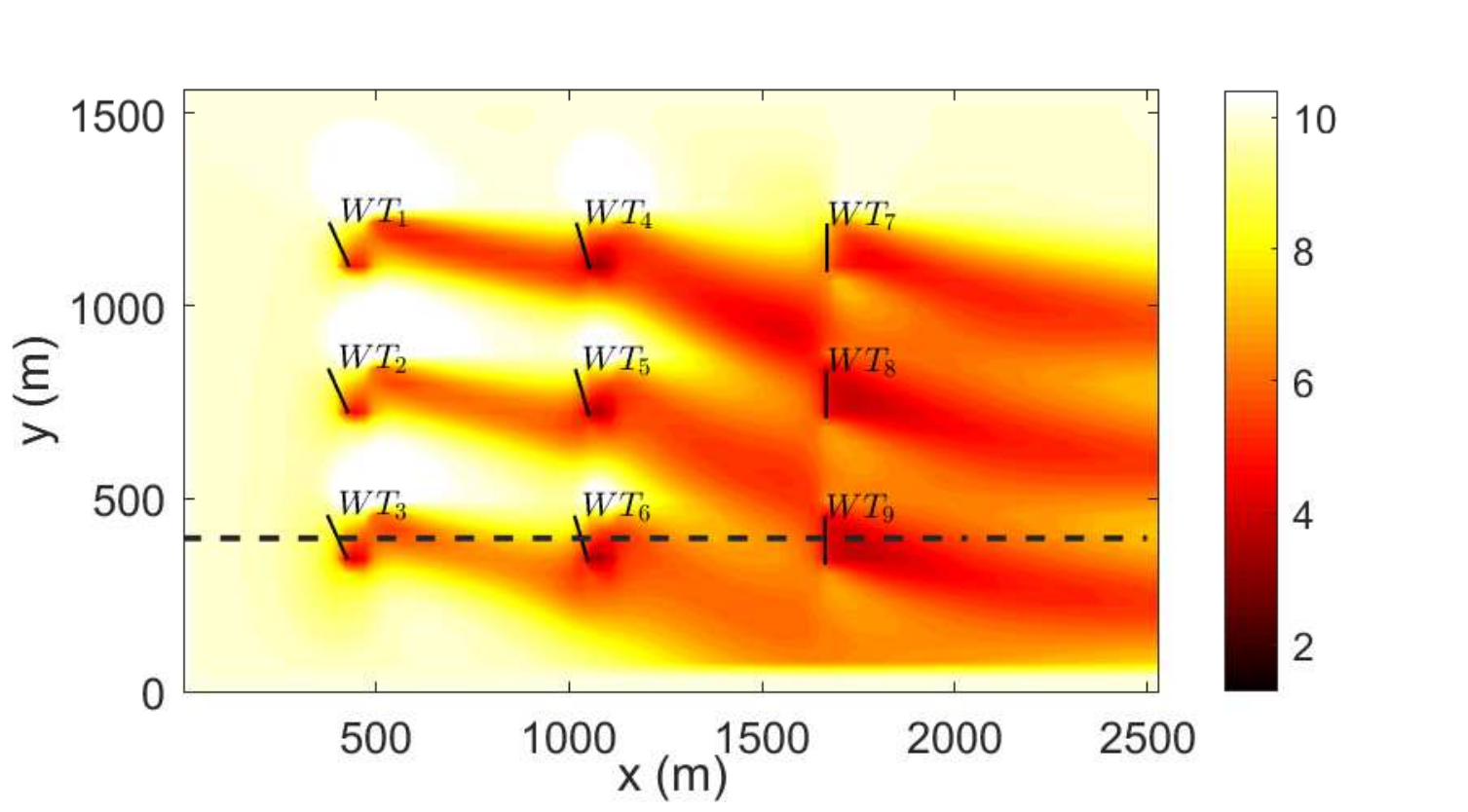}
    \vspace{-7.5pt}
	\caption{Instantaneous streamwise flow velocity snapshot at $t=600s$.}
	\label{fig:controllaw}
\end{figure} 

\begin{figure}[h]
	\centering
		\includegraphics[width=3.5 in]{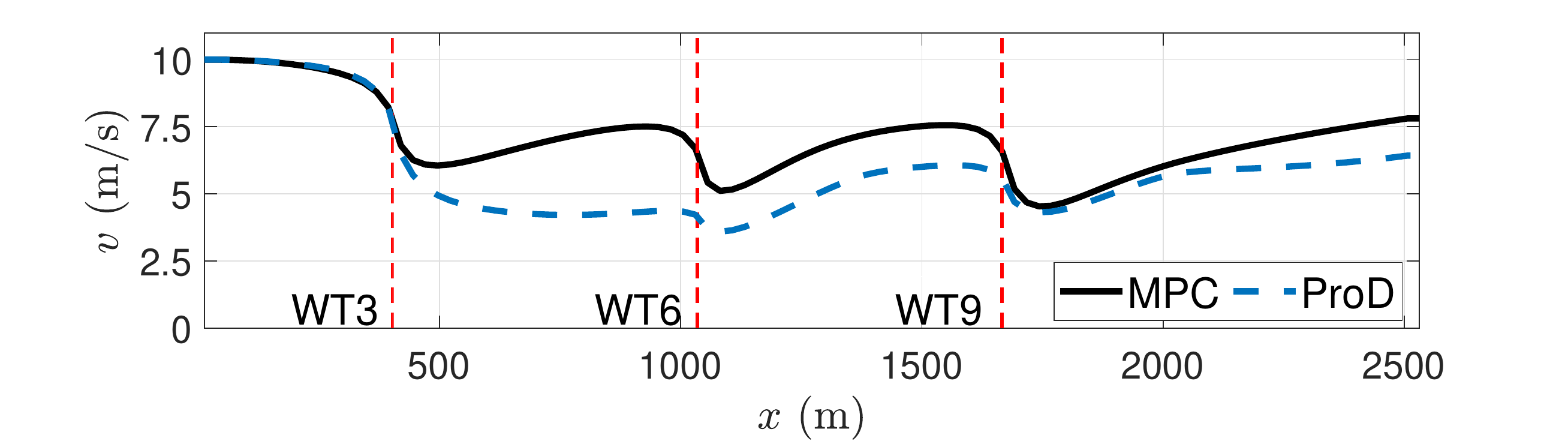}
   \vspace{-7.5pt}
	\caption{Wind speed profile at $t=600s$ for the line of $y=400m$.}
	\label{fig:speedin}
\end{figure}

\begin{figure}[h]
	\centering
		\includegraphics[width=3.25in]{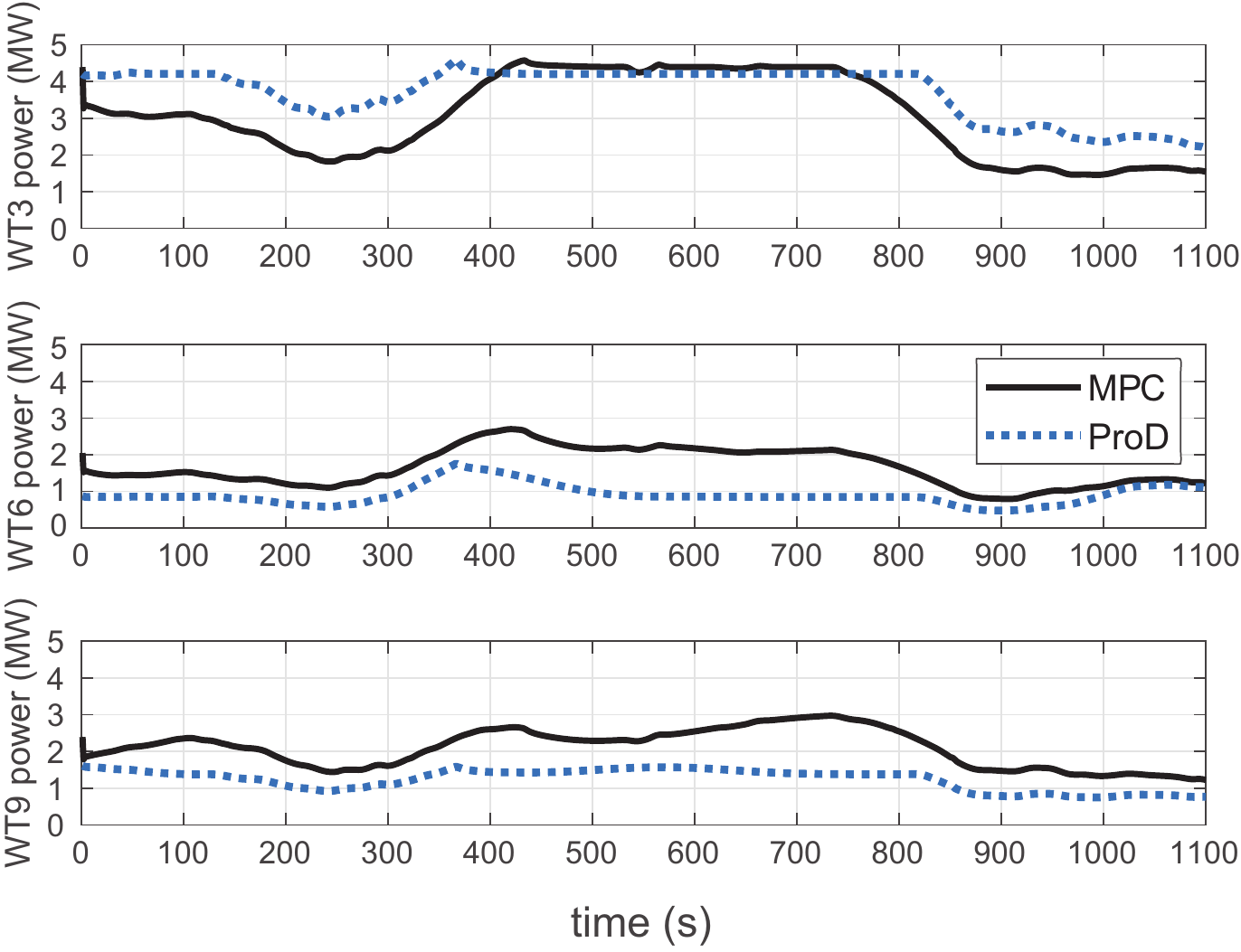}
    \vspace{-7pt}
	\caption{Turbine power production controlled by MPC and ProD.}
    \vspace{-4.5pt}
	\label{fig:indivpower}
\end{figure}

Fig. \ref{fig:controllaw} and \ref{fig:indivpower} further depict the coordinated operation of different turbines for the Fig. \ref{fig:AGCgreedy} (b) case. The yaw angles from upstream to downstream turbines gradually swing to $25^\textrm{o}$, $16.3^\textrm{o}$ and $0^\textrm{o}$ by the proposed MPC. As shown in Fig. \ref{fig:controllaw}, the upstream turbines purposely yaw to deflect the wake so that it will not at all or partially overlap the downwind turbines and more power can be harvested. The detailed wind speed profile is shown in Fig. \ref{fig:speedin}. The instantaneous wind speed at $t=600s$ on the dashed line of $y=400\textrm{m}$ in Fig. \ref{fig:controllaw} is plotted. By coordinating WF operation, the wind speed faced by downstream turbines in the MPC case is higher compared to the ProD case.
Fig. \ref{fig:indivpower} additionally depicts the individual turbine power. The power production of individual WTs is adjusted dynamically by the WF controller to track the power reference. In both cases, the upstream turbine produces the most power since the wind speed in front of these turbines is the highest. Furthermore, it is proven that more power can be harvested by the downstream WTs in the MPC case. 

Through the proposed WF ROM and MPC framework, the trackable power signal range increases to the maximal power extraction of the WF when it operates at its global optimal point considering wake interactions. For the considered 9-turbine WF, the range of trackable AGC signals is improved by $50\%$ compared with traditional ProD control. 

With the reduction of system state dimensions, the WF ROM is computational efficient for control applications.
On a 3.0 GHz Intel Core with i9 processor, for the 9-turbine WF case, the averaged computation time per simulation time step of WFSim is 167 ms. On the other hand, for the WF ROM with 20 states, the computation time for one time step is only 0.452 ms, which is around $0.3\%$ of the time taken by WFSim. For optimally controlling the WF ROM with the prediction horizon $T_p = 250s$, each optimization iteration in Algorithm 1 takes approximately 0.48 s with the CVXOPT toolbox in Python. This speed is orders of magnitude faster than optimal control directly applied to the full-order system, witch takes approximately 133 s per iteration \cite{vali2019adjoint}.

\subsection{25-turbine WF Example} \label{25wtcase}
\begin{figure}[h]
	\centering
		\includegraphics[width=3.3 in]{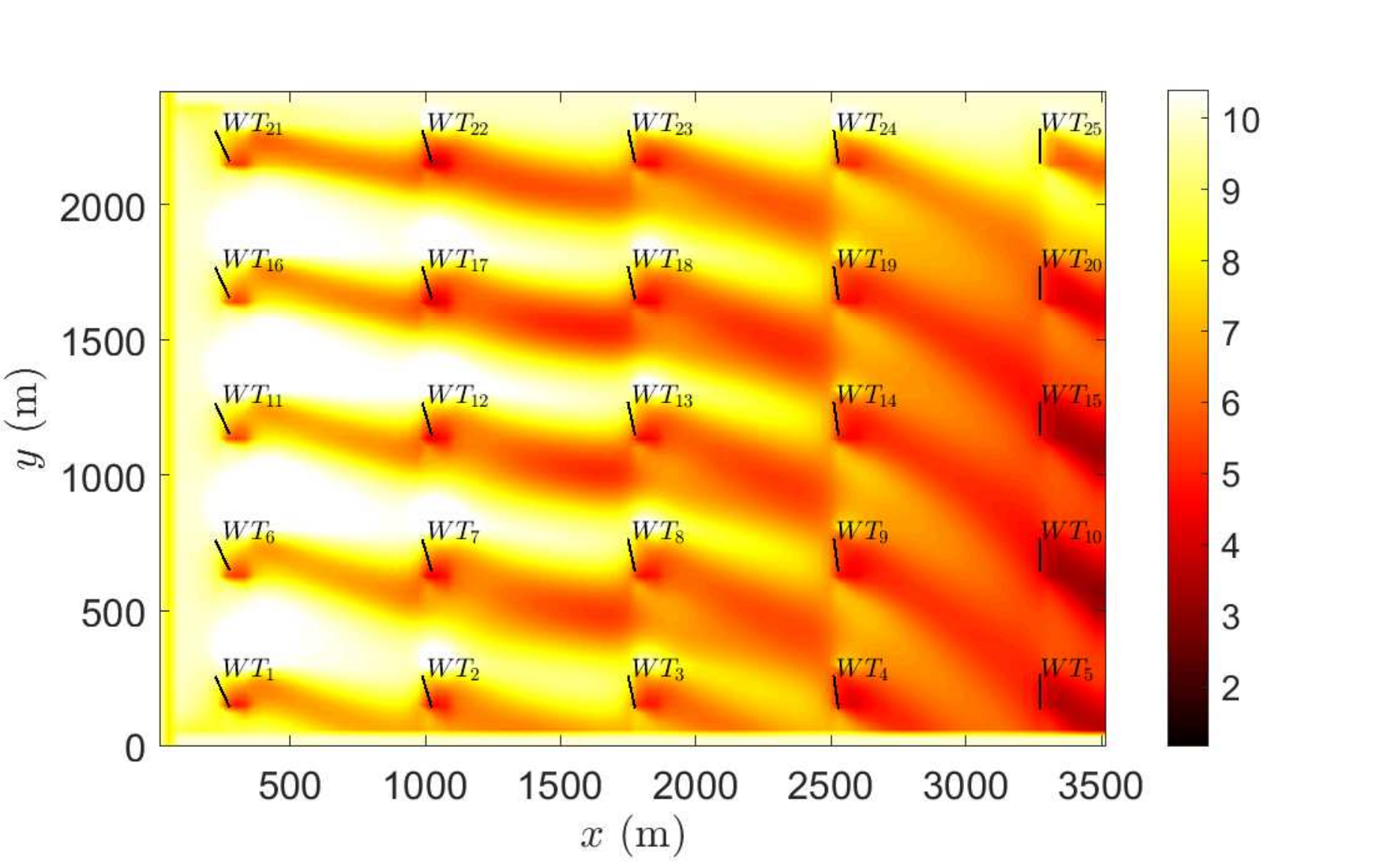}
	\caption{Instantaneous streamwise velocity at $t=600s$ for 25-turbine WF}
	\label{fig:controllaw_25WT}
\end{figure}

\begin{figure}[h]
	\centering
	\includegraphics[width=3.3 in]{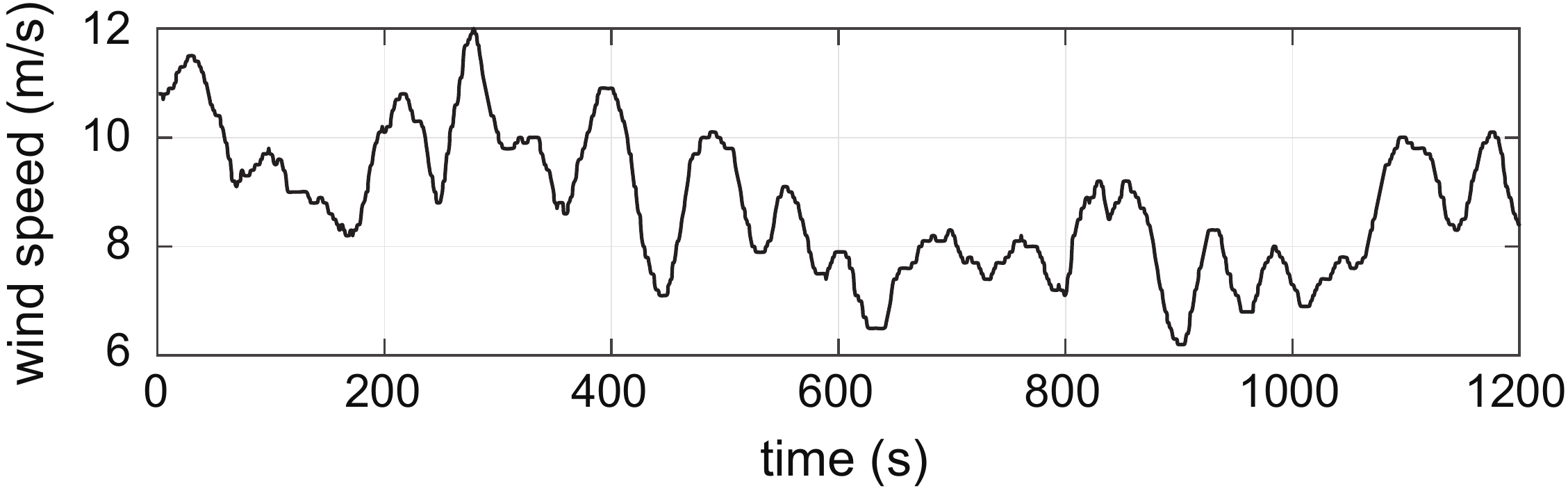}
    \vspace{-8pt}
    \caption{Varying inflow wind speed.}
	\label{fig:speed_25WT}
\end{figure}

A larger WF is further studied with 25 turbines that have 6D spacing in the streamwise direction and 4D in the spanwise direction, as shown in Fig. \ref{fig:controllaw_25WT}. To our best knowledge, this is the largest WF case in WF control studies considering the dynamic wake effects. For even larger wind farms, wind farm partitioning methods can be exploited to split the wind farm into smaller parts \cite{SINISCALCHIMINNA2020656}, and the proposed methods can be applied on individual ones.  
Furthermore, the performance of the deep learning-aided MPC is evaluated with a varying incoming wind speed. The incoming wind speed is shown in Fig. \ref{fig:speed_25WT}, which is a period of realistic measurement data in Denmark and is applied to the west boundary of $x=0\textrm{m}$ as inflow boundary condition. The file of case settings and the wind speed measurement data are also provided online\footnote{The data has been uploaded to https://github.com/kkxchen/WF-deep-learning-aided-MPC-supplementary-material for easy access} .

In this simulation, the simulated wind field is $3500 \times 2416 \textrm{m}^2$, with a grid of $100\times82$ cells ($N_x \times N_y$). The corresponding cell size of the discretization is $35m \times 30m$. The position of first turbine at the left bottom is (250m,200m).
A flow field is first generated with the uniform wind speeds of 10 m/s in the longitudinal direction. Then, for the considered topology, the flow is propagated 400 seconds in advance with the greedy control of all turbines so that the wakes are fully developed. The flow field obtained after the initialization is utilized as the initial flow field for the simulation results presented in this work. In this experiment, the network architecture of the encoder is shown in Table \ref{Tab2}.  
\begin{table}[t]
\centering
\caption{Network Structure of the Encoder in the 25-WT Case}
\label{Tab2}
\begin{tabular}{c c c }
\hline
  Block name & Output size & Layer parameters \\
\hline
  Res Block.1 & $50 \times 41$ & \fontsize{9pt}{10pt} $\left[\begin{smallmatrix}  1 \times 1, 64 \\  3 \times 3, 64 \\  1 \times 1, 256 \end{smallmatrix}\right] \times 3$ \\
  Res Block.2 & $25 \times 21$ &\fontsize{9pt}{10pt}  $\left[\begin{smallmatrix} 1 \times 1, 128 \\  3 \times 3, 128 \\  1 \times 1, 512 \end{smallmatrix}\right] \times 3$ \\
  Res Block.3 & $13 \times 11$ &\fontsize{9pt}{10pt}  $\left[\begin{smallmatrix} 1 \times 1, 256 \\  3 \times 3, 256 \\  1 \times 1, 1024 \end{smallmatrix}\right] \times 3$ \\
  Dense & $n_z=50$ & $\begin{smallmatrix} 4\hspace{2pt}\textrm{fully} \hspace{2pt}\textrm{ connected} \hspace{2pt} \textrm{layers} \\ 4000, 2000, 1200, 500 \end{smallmatrix}$ \\
\hline
\end{tabular}
\end{table}

\begin{figure}[h]
	\centering
		\includegraphics[width=3.3 in]{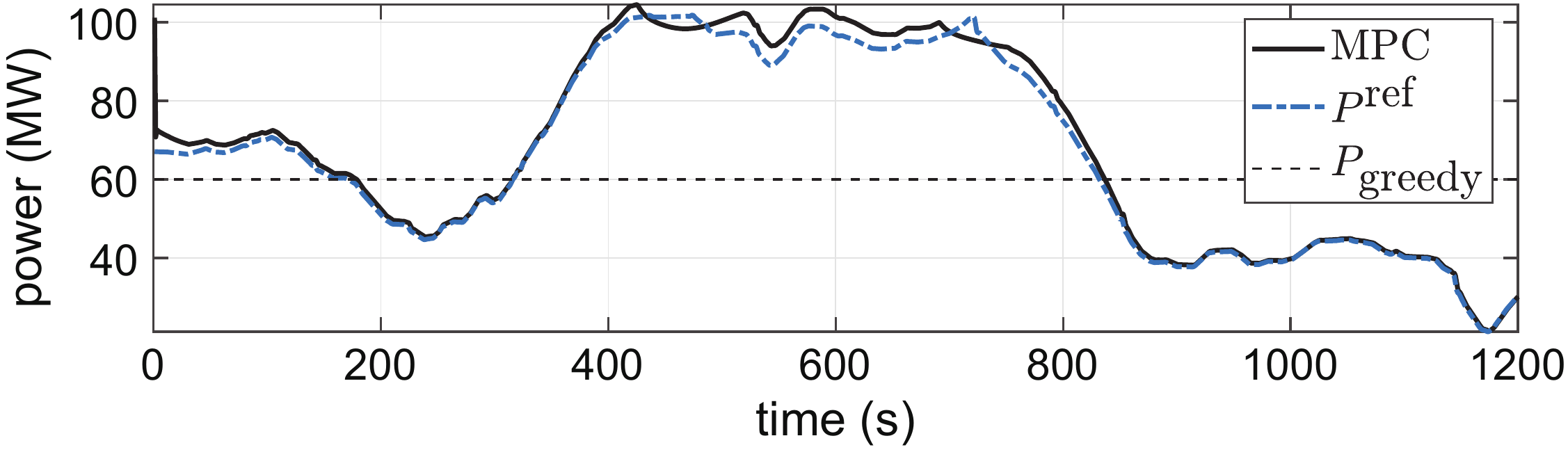}
	\caption{WF power tracking performance for the 25-turbine WF.}
	\label{fig:AGC_25WT}
\end{figure}

The control performance of the proposed method under the varying inflow case is shown in Fig. \ref{fig:AGC_25WT}. Through the proposed framework, the complex flow dynamics can be considered in the optimal control efficiently. Compared to the greedy operation, the trackable power signal range is improved by $68\%$ in this case. Besides, the feedback mechanism of MPC makes the proposed control method robust to the incoming wind variation. At every $T_a$ time, the current wind field full states are observed and fed back to the controller. The ROM prediction error for the wind field introduced by incoming speed turbulence is thus corrected. 

For the 25-turbine WF case, the averaged computation time per simulation time step of WFSim is 210 ms. For the derived WF ROM with 50 states, the computation time per time step is only 1.21 ms. To optimally control the WF ROM with the proposed deep learning-aided MPC, the averaged computational cost is only 1.3 s per iteration in Algorithm 1.

\section{Conclusion}\label{cons}
In this paper, a novel control-oriented WF reduced-order dynamic model is constructed by two specially designed DNN architectures. The original input-output relation is well preserved with good approximation accuracy in the whole control feasible region, while the number of system states is reduced by many orders of magnitude. Besides, a novel deep learning-aided MPC algorithm is proposed to optimally control the learned WF ROM, whose computational advantages are analyzed both theoretically and experimentally. The range of trackable dynamic power production increases to the WF maximal power extraction considering wake interactions. Under the wake-loss-heavy situation, the possible power range that can be tracked in the dynamic operation is improved by $50\%$ in a 9-turbine WF compared with traditional ProD control, agreeing well with the studies in the literature \cite{boersma2017tutorial, munters2018dynamic}. A larger WF with 25 turbines is also studied, showing that the proposed method can obtain similar improvements even in the condition when the incoming wind speed is disturbed.

Future work may include applying this novel method in more complex control objectives like reducing wind turbine loads. Another direction is to incorporate
the three-dimensional wake dynamics in the proposed models and control algorithm.

\section*{References}

  \bibliographystyle{elsarticle-num}
  \bibliography{DLMPC-WF}





\nolinenumbers
\end{document}